\documentclass[11pt]{article}
\usepackage{amsmath, amsfonts, amsthm,amssymb, color, hyperref, enumerate, extarrows, cite}
\usepackage{comment}

\newtheorem{theorem}{Theorem}
\newtheorem{lem}{Lemma}

\newtheorem{definition}[lem]{Definition}
\newtheorem*{notation*}{Notation}
{\theoremstyle{definition}
\newtheorem{example}[lem]{Example}}

\numberwithin{equation}{section}
{\theoremstyle{definition}
\newtheorem{remark}{Remark}}
\usepackage[T1]{fontenc}

\hypersetup{
	colorlinks   = true,
	citecolor    = magenta}

\usepackage[titletoc]{appendix}

\begin{document}

\title{\vspace{-1.2cm} \bf 
Rotational symmetries of domains and orthogonality relations \rm}

\author{Soumya Ganguly, John N. Treuer}
\date{}

\maketitle

\begin{abstract}
Let $\Omega \subset \mathbb{C}^n$ be a domain whose Bergman space contains all holomorphic monomials.  We derive sufficient conditions for $\Omega$ to be Reinhardt, complete Reinhardt, circular or Hartogs in terms of the orthogonality relations of the monomials with respect to their $L^2$-inner products and their $L^2$-norms. More generally, we give sufficient conditions for $\Omega$ to be invariant under a linear group action of an $r$-dimensional torus, where $r \in \{1,\ldots, n\}$. 

\end{abstract}

\renewcommand{\thefootnote}{\fnsymbol{footnote}}
\footnotetext{\hspace*{-7mm} 
\begin{tabular}{@{}r@{}p{16.5cm}@{}}
& Keywords: Bergman space, Moment problems, Reinhardt domains, Circular domains, Quasi-Reinhardt domains.\\
& Mathematics Subject Classification: Primary 32Q06; Secondary 44A60, 32A36. 

\end{tabular}

\noindent\thanks{The second author is supported in part by the NSF grant DMS-2247175 Subaward M2401689. }

\noindent \date
}

\section{Introduction}\label{intro}
Let $\Omega$ be a domain in ${\mathbb C}^n$. Can the rotational symmetries of $\Omega$ be detected from the orthogonality relations of the monomials in the Bergman space? We answer this question for a broad class of domains that include Reinhardt, circular, and Hartogs domains (see Definition \ref{Reinhardt definition} below) satisfying an additional restriction on the $L^2$-norms of the monomials. The key tool is a standard result on the complex moment problem. 

Let $A^2(\Omega)$ denote the Bergman space of the $L^2$-holomorphic functions on $\Omega$. The Bergman space is a separable Hilbert space with a complex inner product 
\begin{equation}\label{A2 inner product}
\langle f, g\rangle := \int_{\Omega} f\bar{g}\, dv, \quad f, g \in A^2(\Omega),
\end{equation}
where $dv$ denotes the Lebesgue volume measure on $\mathbb{C}^n = \mathbb{R}^{2n}$. The symmetry relations of the Reinhardt, circular and Hartogs domains are reflected in the orthogonality relationships of the monomials in the Bergman space with respect to the inner product \eqref{A2 inner product}: 

\begin{enumerate}[i.]
     \item \label{statement 1} If $\Omega \subset \mathbb{C}^n$ is Reinhardt and $z^{\alpha}, z^{\beta} \in A^2(\Omega)$ with $\alpha \neq \beta$, then $z^{\alpha}$ and $z^{\beta}$ are orthogonal.
     \item \label{statement 2} If $\Omega \subset \mathbb{C}^n$ is circular and $z^{\alpha}, z^{\beta} \in A^2(\Omega)$ satisfy that $|\alpha| \neq |\beta|$, then $z^{\alpha}$ and $z^{\beta}$ are orthogonal. 
     \item \label{statement 3} If $\Omega \subset \mathbb{C}^n$ is Hartogs in the $j^{th}$-coordinate and $z^{\alpha}, z^{\beta} \in A^2(\Omega)$ with $\alpha_j \neq \beta_j$, then $z^{\alpha}$ and $z^{\beta}$ are orthogonal. 
 \end{enumerate}
As standard, in statements \ref{statement 1}-\ref{statement 3}, $\alpha = (\alpha_1, \ldots, \alpha_n) \in \mathbb{N}^n$, where $\mathbb{N}$ is the set of all natural numbers including 0, $z^{\alpha} = z_1^{\alpha_1}\cdots z_n^{\alpha_n}$, and $|\alpha| = \alpha_1 + \cdots +\alpha_n$. Similar conventions apply to the multi-index $\beta$. 

Proofs of these facts will naturally follow from equation \eqref{calculation} below. In this paper, we will consider converses of statements \ref{statement 1}.-\ref{statement 3}.~for a much broader class of domains with rotational symmetries called $\rho_A$-invariant domains. These domains, originally considered by Deng and Rong \cite{DeRo16} and Li and Rong \cite{LiRo19}, generalize the notions of Reinhardt, circular, and Hartogs domains.

\begin{definition}\label{rhoA notation definition}
    Let $S = \{\zeta \in \mathbb{C}: |\zeta| = 1\}$.  For $r \in \{1,\ldots, n\}$, let $A \in M_{n \times r}(\mathbb{Z})$ and $\lambda = (\lambda_1,\ldots, \lambda_r) \in S^r$.  Define the group action (of the $r$-dimensional torus) $\rho_A$ on $\mathbb{C}^n$ by
\begin{equation}\label{definition of rho_A}
 \rho_A(\lambda)z = (\lambda_1^{a_{11}}\lambda_2^{a_{12}}\cdots\lambda_r^{a_{1r}}z_1,\ldots, \lambda_1^{a_{n1}}\cdots\lambda_r^{a_{nr}}z_n) = \hbox{Diag}(\lambda_1^{a_{11}}\lambda_2^{a_{12}}\cdots\lambda_r^{a_{1r}},\ldots, \lambda_1^{a_{n1}}\cdots \lambda_r^{a_{nr}})\cdot z \end{equation}
 where $\cdot$ denotes matrix multiplication and $z \in \mathbb{C}^n$ is interpreted as a column vector. A domain $\Omega$ is $\rho_A$-invariant if $\rho_A(\lambda)\Omega := \{\rho_A(\lambda)z : z \in \Omega\} = \Omega$ for every $\lambda \in S^r$.
\end{definition}

We formally define Reinhardt, (quasi-)circular, and Hartogs domains in terms of a group action of the torus.

 \begin{definition}\label{Reinhardt definition}
      Using the notation of Definition \ref{rhoA notation definition},
      \begin{enumerate}
      
     \item  when $r = n$ and $A$ is the $n \times n$ identity matrix, 
     $$
     \rho_A(\lambda)(z) = (\lambda_1z_1,\ldots, \lambda_nz_n).
     $$
     In this case, the $\rho_A$-invariant domains are called Reinhardt domains. A complete Reinhardt domain is a Reinhardt domain that is also star-shaped with respect to the origin. 
     \item  when $r = 1$, i.e.~$\lambda \in S^1$, and $A = (m_1, \ldots, m_n)^T$, with $gcd(m_1, \ldots, m_n)=1$, and $m_i > 0$ for all $i$,
     $$
     \rho_A(\lambda)(z) = (\lambda^{m_1}z_1,\ldots, \lambda^{m_n}z_n).
     $$
     In this case, the $\rho_A$-invariant domains are called quasi-circular domains with weight $(m_1,\ldots, m_n)$.  In particular, if $m_1 = \cdots = m_n = 1$, the $\rho_A$-invariant domains are called circular domains.   
     \item when $r = 1$, i.e.~$\lambda \in S^1$, and $A = (m_1, \ldots, m_n)^T$, with $m_k$ = 0 for $k \neq j$ and $m_j = 1$, 
     \begin{equation}\label{rho for Hartogs domains}
     \rho_A(\lambda)(z) = (z_1,\ldots, z_{j-1}, \lambda z_j, z_{j + 1},\ldots, z_n).
     \end{equation}
     In this case, the $\rho_A$-invariant domains are called Hartogs domains in the $j^{th}$-coordinate.
 \end{enumerate}
\end{definition}

\begin{remark}
    When $\Omega$ is $\rho_A$-invariant and the only entire functions on $\mathbb{C}^n$ that are $\rho_A$-invariant are the constant functions, $\Omega$ is called \emph{quasi-Reinhardt} (see \cite{DeRo16} and \cite{LiRo19}). The Reinhardt and quasi-circular domains are quasi-Reinhardt. However, in $\mathbb{C}^n$ for $n > 1$, the Hartogs domains in the $j^{th}$-coordinate are not. Indeed, for $k \neq j$, the function $f(z) = z_k$ is $\rho_A$-invariant as seen in \eqref{rho for Hartogs domains} but not constant. 
\end{remark}

\begin{remark}
    The Bergman spaces of Reinhardt domains have been studied in the context of the Bergman kernel, the domain function $K:\Omega \times \Omega \to \mathbb{C}$ defined by 
    $ K(z, w) = \sum_{j = 1}^{\infty}\phi_j(z)\overline{\phi_j(w)},$
    where $\{\phi_j\}_{j = 1}^{\infty}$ is an orthonormal basis for $A^2(\Omega)$.  Since any Reinhardt domain admits an orthonormal basis of $A^2(\Omega)$ consisting of Laurent monomials, the series expansion of the Bergman kernel has been effectively used in many instances to give explicit formulations of the Bergman kernel on Reinhardt domains. The readers are referred to \cite{DAngelocomplexellipse, FrancsicsHanges1996 } and the references therein for a few examples. See also \cite{Beberok2, park1, debraj1} and their references for more recent works and related techniques for computing the Bergman kernel.
\end{remark}

 For a fixed $A \in M_{n \times r}(\mathbb{Z})$ and each $\alpha, \beta \in \mathbb{N}^n$, set $g_{\alpha, \beta}:S^r \to \mathbb{C}$ to be
 \begin{equation}\label{definition of g}
 g_{\alpha, \beta}(\lambda) =g_{\alpha, \beta}(\lambda_1,\ldots, \lambda_r) := \prod_{k=1}^r\lambda_k^{\sum_{j=1}^na_{jk}(\alpha_j - \beta_j)}.
 \end{equation}
 When $\Omega \subset \mathbb{C}^n$ is $\rho_A$-invariant,
 \begin{eqnarray}
     \langle z^{\alpha}, z^{\beta}\rangle = \int_{\Omega} z^{\alpha}\bar{z}^{\beta}dv(z) &=& \int_{\Omega}\lambda_1^{a_{11}(\alpha_1 - \beta_1) + \cdots + a_{n1}(\alpha_n - \beta_n)}\cdots\lambda_r^{a_{1r}(\alpha_1 - \beta_1)+\cdots + a_{nr}(\alpha_n - \beta_n)}z^{\alpha}\bar{z}^{\beta}dv(z) \nonumber
     \\
     &=& g_{\alpha, \beta}(\lambda)\int_{\Omega} z^{\alpha}\bar{z}^{\beta} dv(z) = g_{\alpha, \beta}(\lambda) \langle z^{\alpha}, z^{\beta}\rangle. \label{calculation}
 \end{eqnarray}
 This tells us that when $g_{\alpha, \beta}(\cdot)$ is not identically equal to 1, $\langle z^{\alpha}, z^{\beta} \rangle = 0$.  In particular, specializing to the cases of Reinhardt, circular and Hartogs domains implies statements \ref{statement 1}.-\ref{statement 3} above. 

 Let $\|\cdot\|_{L^2(\Omega)} = \|\cdot\|_{L^2(\mathbb{C}^n, \chi_{\Omega}dv)}$ denote the $L^2$-norm on $\Omega$  or equivalently the $L^2$-norm on $\mathbb{C}^n$ with respect to the measure $\chi_{\Omega}dv$.  Here $\chi_{\Omega}$ is the indicator function of $\Omega$. We will prove the converses of statements \ref{statement 1}.-\ref{statement 3}.~for domains that are equal to the interior of their closures and satisfy the following extra assumption.
 
\begin{definition}\label{Condition D}
    Let $\Omega \subset \mathbb{C}^n$ be a domain such that $\{z^{\alpha}\}_{\alpha \in \mathbb{N}^n} \subset A^2(\Omega)$.  For $j \in \{1,\ldots, n\}$, $\Omega$ satisfies Condition $D_j$ if the following series of $L^2$-norms diverges:
\begin{equation}\label{Condition D series which diverges}
\sum_{k = 1}^\infty \|z_j^k\|_{L^2(\Omega)}^{-{1 \over k}} = \infty.
\end{equation}
  We say $\Omega$ satisfies  Condition $D$ if it satisfies Condition $D_j$ for all $j = 1,\ldots, n$.  
\end{definition} 


\begin{remark}\label{bounded domains satisfy Condition D}
    Let $\pi_j:\mathbb{C}^n \to \mathbb{C}$ denote the projection $\pi_j(z_1, \ldots, z_n) = z_j$.  It is easy to see that whenever $\pi_j(\Omega)$ is bounded, Condition $D_j$ holds. In particular, all bounded domains satisfy Condition $D$. In Section \ref{sec:The complex moment problem}, we construct two unbounded, complete Reinhardt domains of the form $\{|z_2| < f(|z_1|)\}$ where $f$ has exponential-type decay, one of which satisfies Condition $D$ and the other which does not but has a Bergman space that contains all holomorphic monomials.
\end{remark}

\begin{remark}
    While a significant portion of the literature on Bergman spaces pertains to those of bounded domains, unbounded domains not biholomorphically equivalent to any bounded one have also been of interest. As a first example, in \cite[Corollary~1.2]{HL24}, Huang and Li constructed an unbounded Reinhardt domain, and determined its Bergman space, whose Bergman metric's holomorphic sectional curvature is identically equal to two.  They also proved that no bounded domain's Bergman metric could have constant holomorphic sectional curvature equal to a positive constant, which demonstrated how the unbounded domain case can differ significantly from the bounded domain case.  As a second example, in \cite{Wi84}, Wiegerinck constructed an unbounded non-pseudoconvex Reinhardt domain whose Bergman space is finite-dimensional.  Such a domain can never be biholomorphic to a bounded domain, and it remains open whether there exists an unbounded pseudoconvex domain with a finite-dimensional Bergman space.

    
\end{remark}
We now state our main results. First, we present our converses. 

\begin{theorem}\label{main Reinhardt theorem}
Let $\Omega \subset \mathbb{C}^n$ be a domain satisfying Condition $D$. 
\begin{enumerate}
    \item \label{case 1} If $\{z^{\alpha}\}_{\alpha \in \mathbb{N}^n}$ is an orthogonal subset of $A^2(\Omega)$ (i.e.~orthogonal with respect to the inner product \eqref{A2 inner product}), then $\hbox{int}(\overline{\Omega})$, the interior of the closure of ${\Omega}$, is Reinhardt. 
    
    \item \label{case 2} If any two monomials $z^{\alpha}$ and $z^{\beta}$ with $|\alpha| \neq |\beta|$ are orthogonal in $A^2(\Omega)$, then $\hbox{int}(\overline{\Omega})$ is circular.
    
    \item \label{case 3}If for some $j \in \{1,\ldots,n\}$, the monomials $z^{\alpha}$ and $z^{\beta}$ are orthogonal in $A^2(\Omega)$ whenever $\alpha_j \neq \beta_j$, then $\hbox{int}(\overline{\Omega})$ is Hartogs in the $j^{th}$ coordinate.
\end{enumerate}
In particular in cases \ref{case 1}.-\ref{case 3}., if $\Omega=\hbox{int}(\overline{\Omega})$, then $\Omega$ is Reinhardt, circular or Hartogs in the $j^{th}$-coordinate respectively.
\end{theorem}
\begin{remark}
    The conclusions that $\hbox{int}(\overline{\Omega})$ is respectively Reinhardt, circular or Hartogs in the $j^{th}$-coordinate cannot be weakened in general to $\Omega$ satisfying those respective conditions. For example, let $\Omega$ be a ball centered at the origin less a point $p=(p_1, \ldots, p_n)$ such that $p_j \ne 0$. Then $\Omega$ satisfies the hypotheses of the three cases above but is neither Reinhardt, circular nor Hartogs in the $j^{th}$-coordinate.
\end{remark}

In case \ref{case 1} of Theorem \ref{main Reinhardt theorem} above, when $\Omega = \hbox{int}(\overline{\Omega})$, $\Omega$ is additionally a domain of holomorphy and the monomials are dense in the Bergman space, we can conclude that $\Omega$ is complete Reinhardt.

\begin{theorem}\label{Complete Reinhardt theorem}
Let $\Omega \subset \mathbb{C}^n$ be a domain satisfying Condition $D$ and suppose that $\{z^{\alpha}\}_{\alpha \in \mathbb{N}^n}$ is an orthogonal basis of $A^2(\Omega)$.  If $\Omega = \hbox{int}(\overline{\Omega})$ and $\Omega$ is a domain of holomorphy, then $\Omega$ is complete Reinhardt.
\end{theorem}

Theorem \ref{main Reinhardt theorem} above, in fact follows from a more general result:

\begin{theorem}\label{main theorem second formulation}
     Suppose $\Omega \subset \mathbb{C}^n$ is a domain satisfying Condition $D$.  Let $r \in \{1,\ldots, n\}$, $A \in M_{n \times r}(\mathbb{Z})$ and for each $\alpha$ and $\beta$, define $g_{\alpha, \beta}$ as in \eqref{definition of g}. Suppose that $g_{\alpha, \beta}(\cdot)$ is identically equal to 1 when $\langle z^{\alpha}, z^{\beta} \rangle \neq 0$.  Then $\hbox{int}(\overline{\Omega})$ is $\rho_A$-invariant.
 \end{theorem}


Theorem \ref{main theorem second formulation} is formulated in terms of the action $\rho_A$, because it can be specialized to many different types of rotational symmetries. For example in Section 3 we use Theorem \ref{main theorem second formulation} to construct orthogonality relations of monomials that are sufficient for a domain to have both Reinhardt and quasi-circular symmetries. 

The proof of Theorem \ref{main theorem second formulation} utilizes a uniqueness result from the complex moment problem that depends on Condition $D$.  In Section \ref{sec:The complex moment problem}, we give a brief introduction to the complex moment problem, and in particular, explain how Condition $D$ fits into this theory. In Section 3, we prove Theorems \ref{main Reinhardt theorem}, \ref{Complete Reinhardt theorem}, and \ref{main theorem second formulation}. 

\subsection*{Acknowledgement}
     The authors are grateful to the anonymous referee for valuable comments on the manuscript and Ziming Shi for insightful conversations regarding this paper.

\section{The complex moment problem}\label{sec:The complex moment problem}
The complex moment problem asks given $\{s_{\alpha\beta}\}_{\alpha, \beta \in \mathbb{N}^n}$ does there exist a unique nonnegative Borel measure $d\mu$ such that 
\begin{equation}\label{complex moment problem}
\int_{\mathbb{C}^n}z^{\alpha}\bar{z}^{\beta}d\mu(z) = s_{\alpha\beta}, \quad \alpha, \beta \in \mathbb{N}^n?
\end{equation}
The sequence $\{s_{\alpha\beta}\}_{\alpha, \beta \in \mathbb{N}^n}$ is called the moment data for the complex moment problem \eqref{complex moment problem}.  When the measure $d\mu$ exists, a sufficient condition for uniqueness is as follows.  

\begin{theorem}(\hspace{-1 pt}\cite[Theorem 12]{StSz85}, \cite[Theorem 15.11]{Sc17})\label{uniqueness theorem of the complex moment problem}
Suppose $\{s_{\alpha\beta}\}_{\alpha, \beta \in \mathbb{N}^n}$ is a sequence of numbers such that $d\mu$ is a nonnegative Borel measure which solves the complex moment problem \eqref{complex moment problem}.  If
\begin{equation}\label{Carleman-type determinacy condition}
\sum_{k = 1}^\infty \|z_j^k\|_{L^2(\mathbb{C}^n, d\mu)}^{-{1 \over k}} = \infty,
\end{equation}
then $d\mu$ is the unique nonnegative Borel measure solving \eqref{complex moment problem}.  
\end{theorem}

Recently, a particular case of Theorem \ref{uniqueness theorem of the complex moment problem} was used by Huang, Li and the second present author in \cite{HL24} to prove that there is no complex manifold whose Bergman space is base point free, separates holomorphic directions and separates points, and whose Bergman metric has identically zero holomorphic sectional curvature. Notice that when $d\mu = \chi_{\Omega}dv$, \eqref{Carleman-type determinacy condition} is the same as \eqref{Condition D series which diverges} in the definition of Condition $D$.  In this paper, when applying Theorem \ref{uniqueness theorem of the complex moment problem}, we will always have $d\mu = \chi_{\Omega}dv$.  Hence, Theorem \ref{uniqueness theorem of the complex moment problem} can be specialized as `\emph{If $d\mu = \chi_{\Omega}dv$ solves the complex moment problem \eqref{complex moment problem}, and $\Omega$ satisfies Condition $D$, then $d\mu$ is the unique solution to \eqref{complex moment problem}'.} 

As stated in Remark \ref{bounded domains satisfy Condition D}, all bounded domains satisfy Condition $D$; however, not all unbounded domains satisfy Condition $D$. Condition $D$ can be violated in two ways.  Either an unbounded domain's Bergman space does not contain all holomorphic monomials, or it contains them but \eqref{Carleman-type determinacy condition} does not hold for some $j$. 
 Let $\Omega = \{(z_1, z_2) \in \mathbb{C}^2:\, |z_2| < f(|z_1|)\}$ for some positive function $f$.  If the decay rate of $f$ is bounded below by a power function, that is, if there are constants $p, C_1, C_2 > 0$ such that $f(|z_1|) \geq C_1|z_1|^{-p}$ for $|z_1| \geq C_2$, then not all monomials will be in the Bergman space.  Indeed
 $$
 \int_{\Omega} |z_1|^{2j} dv(z_1, z_2) = 2 \pi^2 \int_0^{\infty} |z_1|^{2j+1}f(|z_1|)^2 d|z_1| \geq 2C_1^2\pi^2  \int_{C_2}^{\infty}|z_1|^{2j+1-2p} d|z_1|= \infty,
 $$
 whenever $j \geq p-1$.  On the other hand, when the decay of $f$ is exponential, Condition $D$ may or may not hold depending on the rate of the decay as the following example shows.

\begin{example}
Let $$\Omega_k= \left\{(z_1, z_2) \in \mathbb{C}^2: |z_2| < e^{-|z_1|^{1/2^k}}\right\}, \quad k = 0, 1.$$ These domains are unbounded, complete Reinhardt, and $\{z^{\alpha}\}_{\alpha \in \mathbb{N}^2} \subset A^2(\Omega_k)$, yet only $\Omega_0$ satisfies Condition D. $\Omega_1$ does not satisfy Condition D. 

Indeed one can use the definition of these domains to verify that they are unbounded and complete Reinhardt. Using polar coordinates, we can notice that  
\begin{eqnarray*}
\int_{\Omega_k} |z_1|^j |z_2|^m \, dv(z_1, z_2)
&=& (2\pi)^2 \int_{0}^{\infty} |z_1|^{j + 1}\int_{0}^{e^{-|z_1|^{1/2^k}}} |z_2|^{m+1} d|z_2|d|z_1| 
\\
&=& (2\pi)^2 \int_{0}^{\infty} \frac{|z_1|^{j+1}}{m+2} \left[ e^{-|z_1|^{1/2^k}} \right]^{m+2} d|z_1| 
\\
&=& (2\pi)^2 \int_{0}^{\infty} \frac{|z_1|^{j+1}}{m+2} \left[ e^{-(m+2)|z_1|^{1/2^k}} \right] d|z_1| 
\\
&=& {(2\pi)^2 \over (m + 2)} \int_{0}^{\infty} \left(\frac{u}{m+2}\right)^{2^k(j + 1)}e^{-u} \left(\frac{2^k u^{2^k -1}}{(m+2)^{2^k}}\right) du
\\
&=& \frac{2^{k+2}\pi^2}{(m+2)^{2^k(j+2)+1}} \int_{0}^{\infty} u^{2^k(j+2)-1} e^{-u} du
\\
&=& \frac{2^{k+2}\pi^2}{(m+2)^{2^k(j+2)+1}}(2^k(j+2)-1)!,
\end{eqnarray*}
where the fourth equality used the substitution $u = (m + 2)|z_1|^{1/2^k}$. Setting $j=2 \alpha_1$ and $m=2\alpha_2$, we get $\|z_1^{\alpha_1}z_2^{\alpha_2}\|_{L^2(\Omega_k)} < \infty$ for all $(\alpha_1, \alpha_2) \in \mathbb{N}^2$. The calculation above also tells us that
\begin{align}\label{eqn:z_1momentOmega_k}
    ||z_1^{\alpha_1}||_{L^2(\Omega_k)} = \Bigl[\frac{2^{k+2}\pi^2}{2^{2^k(2\alpha_1+2)+1}}(2^k(2\alpha_1+2)-1)!\Bigr]^{1/2}.
\end{align}
Now we note for both $\Omega_0$ and $\Omega_1$, Condition $D_2$ will be satisfied as $|z_2|<1$ in both domains. However we now show that only $\Omega_0$ satisfies Condition $D_1$. 

When $k=0$, \eqref{eqn:z_1momentOmega_k} gives us  $$
\|z_1^{\alpha_1}\|_{L^2(\Omega_0)}  = \Bigl[{(2\pi)^2(2\alpha_1 + 1)! \over 2^{2\alpha_1 + 3}}\Bigr]^{1/2}. 
$$
Let $C$ denote a numerical constant independent of $\alpha_1$, which may be different at each appearance.  Then 
$$
\sum_{\alpha_1 = 1}^\infty \|z_1^{\alpha_1}\|_{L^2(\Omega_0)}^{-{1 \over \alpha_1}} = \sum_{\alpha_1=1}^{\infty} \left(2^{2\alpha_1 + 3} \over (2\pi)^2(2\alpha_1 + 1)! \right)^{1 \over 2\alpha_1} \geq C\sum_{\alpha_1=1}^{\infty} {1 \over (2\alpha_1 + 1)!^{1 \over 2\alpha_1}} \geq C\sum_{\alpha_1=1}^{\infty}{1 \over (2\alpha_1 + 1)^{2\alpha_1 \over 2\alpha_1}} = \infty,
$$
which shows that Condition $D_1$ is satisfied. So $\Omega_0$ satisfies Condition $D$. 

When $k=1$, \eqref{eqn:z_1momentOmega_k} gives us $$\|z_1^{\alpha_1}\|_{L^2(\Omega_1)}= \Bigl[{\pi^2 (4\alpha_1 + 3)! \over 2^{4\alpha_1 + 2}}\Bigr]^{1/2}.$$ From here we can notice that
\begin{eqnarray}
\sum_{\alpha_1 =1}^{\infty} \|z_1^{\alpha_1}\|_{L^2(\Omega_1)}^{-{1\over \alpha_1}} = \sum_{\alpha_1 = 1}^{\infty} \left( {2^{4\alpha_1 + 2} \over \pi^2(4\alpha_1 + 3)! }\right)^{1 \over 2\alpha_1}
= 4\sum_{\alpha_1 = 1}^{\infty}\left({4 \over \pi^2}\right)^{1 \over 2\alpha_1}\left(1 \over (4\alpha_1 + 3)!\right)^{1 \over 2\alpha_1}. \label{the norm before Stirling's formula is used}
\end{eqnarray}
 By Stirling's formula, for $\alpha_1$ sufficiently large,
\begin{eqnarray}
\left[1 \over (4\alpha_1 + 3)! \right]^{1 \over 2\alpha_1}  &\leq& C\left[ \sqrt{2\pi(4\alpha_1 + 3)}\left({4\alpha_1 + 3 \over e}\right)^{4\alpha_1 + 3} \right]^{-{1 \over 2\alpha_1}} \nonumber
\\
&=& C\left[e^{4\alpha_1 + 3} \over \sqrt{2\pi(4\alpha_1 + 3)}(4\alpha_1 + 3)^{4\alpha_1 + 3} \right]^{1 \over 2\alpha_1} \nonumber
\\
&=& C\left[{e^{2 + {3 \over 2\alpha_1}} \over (2\pi(4\alpha_1 + 3))^{1 \over 4\alpha_1}(4\alpha_1 + 3)^{2 + {3 \over 2\alpha_1}}} \right] \nonumber
\\ &\leq& C{1 \over (4\alpha_1 + 3)^2}. \label{key estimate from Stirlings formula}
\end{eqnarray}
Plugging in \eqref{key estimate from Stirlings formula} into \eqref{the norm before Stirling's formula is used} demonstrates that $\sum_{\alpha_1 =1}^{\infty} \|z_1^{\alpha_1}\|_{L^2(\Omega_1)}^{-{1\over \alpha_1}} < \infty$. Thus $\Omega_1$ does not satisfy Condition D.
\end{example}

\section{Proofs of Theorems \ref{main Reinhardt theorem}, \ref{Complete Reinhardt theorem}, and \ref{main theorem second formulation}}
 \begin{proof}[\textbf{Proof of Theorem \ref{main theorem second formulation}}]
    Define $s_{\alpha\beta} = \langle z^{\alpha}, z^{\beta} \rangle$.  By the change of variables theorem,
     \begin{eqnarray}
         s_{\alpha\beta} = \int_{\mathbb{C}^n}z^{\alpha}\bar{z}^{\beta}\left(\chi_{\Omega}(z)dv(z)\right) &=& \int_{\mathbb{C}^n}g_{\alpha, \beta}(\lambda)z^{\alpha}\bar{z}^{\beta}\left(\chi_{\rho_A(\lambda)^{-1}(\Omega)}(z) dv(z)\right) \nonumber
         \\
         &=& \int_{\mathbb{C}^n} z^{\alpha}\bar{z}^{\beta}\left(\chi_{\rho_A(\bar{\lambda})(\Omega)}(z) dv(z)\right), \quad \lambda \in S^r,
     \end{eqnarray}
     where the final equality follows because $s_{\alpha \beta}$ is  equal to 0 if $g_{\alpha, \beta} \not\equiv 1$. Then each of the measures in $\{\chi_{\rho_A(\bar{\lambda})(\Omega)}(z) dv(z)\}_{\lambda \in S^r}$ solve the same complex moment problem with moment data $\{s_{\alpha\beta}\}$. Since Condition $D$ holds, by Theorem \ref{uniqueness theorem of the complex moment problem},
     $$
     \chi_{\rho_A(\bar{\lambda})(\Omega)}(z) dv(z) = \chi_\Omega(z) dv(z), \quad \lambda \in S^r.
     $$
     Thus,
     $$
     \chi_{\rho_A(\bar{\lambda})(\Omega)}(z) = \chi_\Omega(z), \quad dv\hbox{-almost everywhere}, \quad \lambda \in S^r.
     $$
      Let $\hbox{Ext}(\Omega)$ denote the exterior of $\Omega$.  When $B \subset \hbox{Ext}(\Omega)$ is a ball,
    $$
   \int_B \chi_{\rho_A(\lambda)(\Omega)}dv = \int_B\chi_{\Omega}dv= 0,
    $$
    which implies that $\chi_{\rho_A(\lambda)(\Omega)}|_{B} = 0$ almost everywhere.  Since $\rho_A(\lambda)\Omega$ is open, $\chi_{\rho_A(\lambda)(\Omega)}|_B \equiv 0$, giving us $B \subset \hbox{Ext}(\rho_A(\lambda)(\Omega))$.  Hence $\hbox{Ext}(\Omega) \subset \hbox{Ext}(\rho_A(\lambda)(\Omega))$. By a similar argument, it follows that
    $
    \hbox{Ext}(\Omega) = \hbox{Ext}(\rho_A(\lambda)(\Omega)).
    $
    Taking the complement and then the interior of both sides, gives
    $$
    \hbox{int}(\overline{\Omega}) = \hbox{int}(\overline{\rho_A(\lambda)\Omega}), \quad \lambda \in S^r.
    $$
    Since $\rho_A(\lambda)$ is a homeomorphism on $\mathbb{C}^n$
    $$
    \hbox{int}(\overline{\Omega}) = \rho_A(\lambda)\hbox{int}(\overline{\Omega}), \quad \lambda \in S^r,
    $$
    which means that $\hbox{int}(\overline{\Omega})$ is $\rho_A(\lambda)$-invariant.
 \end{proof}

Now Theorem \ref{main Reinhardt theorem} follows as special cases of Theorem \ref{main theorem second formulation} as shown below. 

\begin{proof}[\textbf{Proof of Theorem \ref{main Reinhardt theorem}}]
    For case \ref{case 1}, set $r = n$, $A$ to be the $n \times n$ identity matrix.  Then 
    $$
    g_{\alpha, \beta}(\lambda) = \lambda_1^{(\alpha_1 - \beta_1)}\cdots\lambda_n^{(\alpha_n - \beta_n)}, \quad \alpha, \beta \in \mathbb{N}^n, \quad \rho_A(\lambda)z = (\lambda_1z_1,\ldots, \lambda_nz_n), \quad \lambda \in S^n.
    $$
    Since $g_{\alpha, \beta}$ identically equals 1 when $\alpha = \beta$, by Theorem \ref{main theorem second formulation}, $\hbox{int}(\overline{\Omega})$ is $\rho_{A}$-invariant.  Therefore, $\hbox{int}(\overline{\Omega})$ is Reinhardt.

    For case 2, set $r = 1$ and $A$ to be the $n \times 1$ matrix where all entries are equal to 1. Then 
    $$
    g_{\alpha, \beta}(\lambda) = \lambda^{\sum_{j=1}^n(\alpha_j - \beta_j)}, \quad \alpha, \beta \in \mathbb{N}^n, \quad \rho_A(\lambda)z = (\lambda z_1, \ldots, \lambda z_n), \quad \lambda \in S^1.
    $$
    Since $g_{\alpha, \beta}$ identically equals 1 when the degree of $z^{\alpha}, z^{\beta}$ are the same, $\hbox{int}(\overline{\Omega})$ is $\rho_A$ invariant; that is circular.
    
    For case 3, set $r = 1$ as before and $A$ to be the $n \times 1$ matrix with 1 in the $j^{th}$-coordinate and 0 in the other coordinates. Then 
    $$
    g_{\alpha, \beta}(\lambda) = \lambda^{\alpha_j - \beta_j}, \quad \alpha, \beta \in \mathbb{N}^n, \quad \rho_A(\lambda)z = (z_1,\ldots, z_{j-1}, \lambda z_j, z_{j+1}, \ldots, z_n), \quad \lambda \in S^1.
    $$
    Since $g_{\alpha, \beta}$ is identically equal to 1 when $\alpha_j = \beta_j$, $\hbox{int}(\overline{\Omega})$ is $\rho_A$-invariant; that is Hartogs in the $j^{th}$-coordinate. 
    \end{proof}

A domain $D \subset \mathbb{C}^n$ is said to be an $L^2_h$-domain of holomorphy if for any pair of domains $D_0, \tilde{D} \subset \mathbb{C}^n$ with $\emptyset \ne D_0 \subset \tilde{D} \cap D, \tilde{D} \not \subset D$,
there exists a function $f \in A^2(D)$ such that $f|_{D_0}$ is not the restriction of a
function $\tilde{f}$ which is holomorphic on $\tilde{D}$. We will now utilize a result about $L^2_h$-domains of holomorphy from \cite{{Pflugjarnicki_1997}} to establish Theorem \ref{Complete Reinhardt theorem}.

\begin{proof}[\textbf{Proof of Theorem \ref{Complete Reinhardt theorem}}]
   By Theorem \ref{main Reinhardt theorem}, $\hbox{int}(\overline{\Omega}) = \Omega$ is Reinhardt. Suppose $\Omega$ is not complete Reinhardt. Let $\tilde{\Omega}$ be the smallest complete Reinhardt domain containing $\Omega$.  
   Since $\{z^{\alpha}\}_{\alpha \in \mathbb{N}^n}$ is an orthogonal basis of $A^2(\Omega)$, for any $f \in A^2(\Omega)$,
   $$ f (z) = \sum_{\alpha \in \mathbb{N}^n} a_{\alpha}z^{\alpha}, \quad z \in \Omega.$$
The domain of convergence of any power series is always a complete Reinhardt domain \cite[Corollary 1.16]{Ra86}. So $f$ defines a holomorphic function on $\tilde{\Omega}$. Since $f \in A^2(\Omega)$ is arbitrarily chosen, $\Omega$ is not an $L^2_h$-domain of holomorphy. Then by \cite[Theorem 3.21]{JaPf87} or \cite[Proposition 2]{Pflugjarnicki_1997}, $\Omega$ is not a domain of holomorphy.
\end{proof}

\begin{remark}
    In Theorem \ref{Complete Reinhardt theorem}, if $\Omega = \hbox{int}(\overline{\Omega})$ is not a domain of holomorphy, then $\Omega$ may not be a complete Reinhardt domain even if $\{z^{\alpha}\}_{\alpha \in \mathbb{N}^n}$ is complete and orthogonal.  Indeed, let $\Omega = P_2 \setminus \overline{P_1}$ be the difference of two polydisks in $\mathbb{C}^n$, $n > 1$ centered around the origin and of polyradius $(2, \ldots, 2)$ and $(1, \ldots, 1)$ respectively.  Then $\Omega=\hbox{int}(\overline{\Omega})$ is Reinhardt and bounded.  It follows that $\{z^{\alpha}\}_{\alpha \in \mathbb{N}^{n}}$ is an orthogonal subset of the Bergman space of $\Omega$.  We will show that it is also complete in $A^2(\Omega)$.  Suppose that $f \in A^2(\Omega)$ and 
    \begin{equation}\label{f is orthogonal to all monomials}
    \langle f, z^{\alpha}\rangle = 0, \quad \alpha \in \mathbb{N}^{n}.
    \end{equation}
    By the Hartogs' extension theorem, $f$ extends to a holomorphic function on $P_2$.  Since $f \in L^2(\Omega)$ and $f \in C(\overline{P_1}) \cap L^2(\Omega)$, we see that $f \in L^2(P_2)$.  Given that $\{z^{\alpha}\}_{\alpha \in \mathbb{N}^n}$ is complete in $A^2(P_2)$, there is a power series $\sum_\alpha a_\alpha z^{\alpha}$ convergent in $P_2$ such that
    $$
    \lim_{N \to \infty} \int_{P_2} |f(z) - \sum_{|\alpha| \leq N} a_{\alpha}z^{\alpha}|^2 dv(z) = 0.
    $$
    By monotonicity,
    $$
    \lim_{N \to \infty} \int_{\Omega} |f(z) - \sum_{|\alpha| \leq N} a_{\alpha}z^{\alpha}|^2 dv(z) = 0.
    $$
    By \eqref{f is orthogonal to all monomials}, 
    $$
    \lim_{N \to \infty}\int_{\Omega} \big(|f(z)|^2 + |\sum_{|\alpha| \leq N} a_{\alpha}z^{\alpha}|^2 \big) dv(z) = 0,
    $$
    which implies that $f = 0$.
\end{remark}

Although the Reinhardt, circular and Hartogs domains are the classes of domains of primary interest, Theorem \ref{main theorem second formulation} can be used to construct additional orthogonality to symmetry relationships in domains that combine different rotational symmetries.
\begin{example}
    Let $I_1 = \{s_1,\ldots, s_p\}$ and $I_2 = \{t_1,\ldots, t_q\}$ where $I_1$ and $I_2$ are nonempty disjoint sets with $I_1 \cup I_2 \subset \{1,\ldots, n\}$. Let $W = \{w_{s_1},\ldots, w_{s_p}\}$ be a weight set where each $w_{s_j}$ is a positive integer weight corresponding to index $s_j$, and $gcd(w_{s_1},\ldots, w_{s_p}) = 1$. Suppose that $\Omega \subset \mathbb{C}^n$ is a domain satisfying Condition $D$. For simplicity, also assume that $\Omega = \hbox{int}(\overline{\Omega})$.  Furthermore, suppose the monomials of the Bergman space satisfy the orthogonality relations
    $$
    \langle z^{\alpha}, z^{\beta} \rangle = 0 \Leftrightarrow \sum_{j=1}^p w_{s_j}(\alpha_{s_j} - \beta_{s_j}) \neq 0, \hbox{ or there exists $j$ with $1 \leq j \leq q$ such that } \alpha_{t_j} \neq \beta_{t_j}. 
    $$
    Set $r = q + 1$ and 
    $
    A = (a_{jk})_{j, k = 1}^{j = n, k = r}$ where
    $$
    a_{jk} = \begin{cases}
        1 & j = t_s,\ k = s,\, (1 \leq s \leq q) 
        \\
        w_{j} & j \in I_1, k = r
        \\
        0 & else.
    \end{cases}
    $$
    Then 
    $$
    g_{\alpha, \beta}(\lambda) = \lambda_1^{(\alpha_{t_1} - \beta_{t_1})}\cdots \lambda_q^{(\alpha_{t_q} - \beta_{t_q})}\lambda_r^{\sum_{j=1}^pw_{s_j}(\alpha_{s_j} - \beta_{s_j})}, \quad \alpha, \beta \in \mathbb{N}^n, 
    $$
    which is identically equal to one when $\langle z^{\alpha}, z^{\beta}\rangle \neq 0.$  Thus, by Theorem \ref{main theorem second formulation}, when 
    $
    \rho_A(\lambda)z = (z_1^{\prime}, \ldots, z^{\prime}_n)
    $ where 
    $$
    z^{\prime}_j = \begin{cases}
        z_j & j \not\in I_1 \cup I_2
        \\
        \lambda_r^{w_j}z_j & j \in I_1 
        \\
        \lambda_sz_{t_s} & j = t_s \in I_2,\, (1 \leq s \leq q),
    \end{cases}
    $$
    $\Omega$ is $\rho_A$-invariant. That is, $\Omega$ is quasi-circular in the variable $z_{s_1},\ldots, z_{s_p}$ with weights $(w_{s_1},\ldots, w_{s_p})$ and Reinhardt in the variables $z_{t_1},\ldots, z_{t_q}$.
\end{example}

\bibliographystyle{amsplain}

\bibliography{bibliography}

\fontsize{11}{9}\selectfont

\vspace{0.5cm}

\noindent s1gangul@ucsd.edu;

 \vspace{0.2 cm}

\noindent Department of Mathematics, University of California San Diego, La Jolla, CA 92093, USA

\vspace{0.6 cm}

\noindent jtreuer@ucsd.edu

 \vspace{0.2 cm}

\noindent Department of Mathematics, University of California San Diego, La Jolla, CA 92093, USA

\end{document}